\documentclass[12pt]{article}

\newfont{\Bbbfont}{msbm10 at 10pt} \newfont{\sBbbfont}{msbm7 at 10pt}

\newcommand{\Bbb}[1]{\mathchoice{\hbox{\Bbbfont #1}}{\hbox{\Bbbfont #1}}%
  {\hbox{\sBbbfont #1}}{\hbox{\sBbbfont #1}}}

\newcommand{\0}{\mathchoice
                  {\hbox{0\kern-3.4pt\lower-.35pt\hbox{\vbox%
                        {\hrule height7pt width0.4pt}}}\hskip3.4pt }
                  {\hbox{0\kern-3.4pt\lower-.35pt\hbox{\vbox%
                        {\hrule height7pt width0.4pt}}}\hskip3.4pt }
                  {\hbox{$\scriptstyle0$\kern-2.7pt\lower-.1pt\hbox{\vbox%
                        {\hrule height4.4pt width0.3pt}}}\hskip2pt }
                  {\hbox{$\scriptscriptstyle0$\kern-2.2pt\lower-.1pt\hbox{\vbox
                        {\hrule height3pt width0.2pt}}}\hskip2pt }}
\newcommand{\1}{\hbox{1\kern-3pt{\rm I}}\kern.5pt }
\newcommand{\pdd}[2]{\frac{\partial#1}{\partial#2}}
\newcommand{\cyrg}{\Gamma}

\newcommand{\Mat}{\mbox{\rm Mat}}
\newcommand{\circledast}{\begin{picture}(12,12)
                  \put(2,1){$\ast$} \put(4.6,3.8){\circle{8}}
                        \end{picture}}

\begin{document}
\pagenumbering{arabic}

\title{Correspondence principle for idempotent calculus and some
computer applications\footnote{Institut des Hautes Etudes Scientifiques,
IHES/M/95/33 Avril 1995\newline
Also published in: Idempotency / J. Gunawardena (editor), Cambridge
University Press, Cambridge, 1998, p. 420--443.}}
\author{G. L. Litvinov and V. P. Maslov}
\date{}

\def\x {\stackrel {\textstyle \otimes}{,}}


\maketitle


\section{Introduction}

\quad This paper is devoted to heuristic aspects of the so-called idempotent
calculus. There is a correspondence between important, useful and interesting
constructions and results over the field of real (or complex) numbers and
similar constructions and results over idempotent semirings in the spirit of
N.~Bohr's correspondence principle in Quantum Mechanics. Idempotent analogs for
some basic ideas, constructions and results in Functional Analysis and
Mathematical Physics are discussed from this point of view. Thus the
correspondence principle is a powerful heuristic tool to apply unexpected
analogies and ideas borrowed from different areas of Mathematics and
Theoretical Physics.

It is very important that some problems nonlinear in the traditional sense (for
example, the Bellman equation and its generalizations and the Hamilton--Jacobi
equation) turn out to be linear over a suitable semiring; this linearity
considerably simplifies the explicit construction of solutions. In this case we
have a natural analog of the so-called superposition principle in Quantum
Mechanics (see~\cite{1}--\cite{3}).

The theory is well advanced and includes, in particular, new integration
theory, new linear algebra, spectral theory and functional analysis.
Applications include various optimization problems such as multicriteria
decision making, optimization on graphs, discrete optimization with a large
parameter (asymptotic problems), optimal design of computer systems and
computer media, optimal organization of parallel data processing, dynamic
programming, discrete event systems, computer science, discrete mathematics,
mathematical logic and so on. See, for example, \cite{4}--\cite{64}. Let us
indicate some applications of these ideas in mathematical physics and
biophysics~ \cite{65}--\cite{70}.

In this paper the correspondence principle is used to develop an approach to
object-oriented software and hardware design for algorithms of idempotent
calculus and scientific calculations. In particular, there is a regular method
for constructing back-end processors and technical devices intended for an
implementation of basic algorithms of idempotent calculus and mathematics of
semirings.  These hardware facilities increase the speed of data processing.
Moreover this approach is useful for software and hardware design in the
general case of algorithms which are not ``idempotent''~\cite{72}.

The paper contains a brief survey of the subject but our list of references is
not complete. Additional references could be found in \cite{4}--\cite{9},
\cite{11}, \cite{14}, \cite{15}, \cite{17}, \cite{19}--\cite{24},
\cite{27}--\cite{29}, \cite{47}, \cite{53}, \cite{63};
the corresponding lists
of references are not complete too but very useful.

The authors are grateful to I.~Andreeva, B.~Doubrov, M.~Gromov, J.~Gunawardena,
G.~Henkin, V.~Kolokoltsov, G.~Mascari, P.~Del Moral, A.~Rodionov, S.~Samborski,
G.~Shpiz, A.~Tarashchan for discussions and support.

The work was supported by the Russian Fundation for Basic
Research (RFBR), Project 96--01--01544.

\section{Idempotent quantization and\newline
  dequantization}

\quad Let $\Bbb R$ be the field of real numbers, $\Bbb R_+$ the subset of all
non-negative numbers. Consider the following change of variables:
\begin{equation}
  u\mapsto w=h \ln u,
\end{equation}
where $u\in\Bbb R_+$, $h>0$; thus $u=e^{w/h}$, $w\in\Bbb R$.  We have got a
natural map
\begin{equation}
  D_h: \Bbb R_+ \to A= \Bbb R\cup\{-\infty\}
\end{equation}
defined by the formula (2.1). Denote by $\0$ the ``additional'' element
$-\infty$ and by $\1$ the zero element of $A$ (that is $\1=0$); of course
$\0=D_h(0)$ and $\1=D_h(1)$. Denote by $A_h$ the set $A$ equipped with the two
operations $\oplus$ (generalized addition) and $\odot$ (generalized
multiplication) borrowed from the usual addition and multiplication in $\Bbb
R_+$ by the map $D_h$; thus $w_1\odot w_2=w_1+w_2$ and $w_1\oplus w_2=h\ln
(e^{w_1/h}+e^{w_2/h})$. Of course, $D_h(u_1+u_2)=D_h(u_1)\oplus D_h(u_2)$ and
$D_h(u_1u_2)=D_h(u_1)\odot D_h(u_2)$. It is easy to prove that $w_1\oplus
w_2=h\ln (e^{w_1/h}+e^{w_2/h}) \to \max\{ w_1,w_2\}$ as $h\to 0$.

Let us denote by $\Bbb R_{\max}$ the set $A=\Bbb R\cup\{-\infty\}$ equipped
with operations $\oplus=\max$ and $\odot=+$; set $\0=-\infty$, $\1=0$.
Algebraic structures in $\Bbb R_+$ and $A_h$ are isomorphic, so $\Bbb R_{\max}$
is a result of a deformation of the structure in $\Bbb R_+$. There is an
analogy to the quantization procedure, and $h$ is an analog for the Planck
constant. Thus $\Bbb R_+$ (or $\Bbb R$) can be treated as a ``quantum object''
with respect to $\Bbb R_{\max}$ and $\Bbb R_{\max}$ can be treated as a
``classical'' or ``semiclassical'' object and as a result of a
``dequantization'' of this quantum object.

{\samepage Similarly denote by $\Bbb R_{\min}$ the set $\Bbb R\cup\{+\infty\}$
  equipped with operations $\oplus=\min$ and $\odot=+$; in this case
  $\0=+\infty$ and $\1=0$. Of course, the change of variables $u\mapsto w=-h\ln
  u$ generates the corresponding dequantization procedure for this case.}

The set $\Bbb R\cup\{+\infty\}\cup\{-\infty\}$ equipped with the operations
$\oplus=\min$ and $\odot=\max$ can be obtained as a result of a ``second
dequantization'' with respect to $\Bbb R$ (or $\Bbb R_+)$. In this case
$\0=\infty$, $\1=-\infty$ and the dequantization procedure can be applied to
the subset of negative elements of $\Bbb R_{\max}$ and the corresponding change
of variables is $w\mapsto v=h\ln(-w)$.

\section{Semirings}

\quad It is easy to check that for
these constructed operations $\oplus$ and $\odot$
the following basic properties are valid for all elements $a$, $b$, $c$:
\begin{equation}
  (a\oplus b)\oplus c=a\oplus (b\oplus c);\qquad (a\odot b)\odot c=a\odot
  (b\odot c);
\end{equation}
\begin{equation}
  \0\oplus a=a\oplus \0=a;\qquad \1\odot a=a\odot\1=a;
\end{equation}
\begin{equation}
  \0\odot a=a\odot\0=\0;
\end{equation}
\begin{equation}
  a\odot(b\oplus c)=(a\odot b)\oplus(a\odot c);\qquad (b\oplus c)\odot
  a=(b\odot a)\oplus (c\odot a);
\end{equation}
\begin{equation}
  a\oplus b=b\oplus a;
\end{equation}
\begin{equation}
  a\oplus a=a;
\end{equation}
\begin{equation}
  a\odot b=b\odot a.
\end{equation}

A set $A$ equipped with binary operations $\oplus$ and $\odot$ and having
distinguished elements $\0$ and $\1$ is called a {\it semiring}, if the
properties (axioms) (3.1), (3.2), (3.3), (3.4) and (3.5) are fulfilled. We
shall suppose that $\0\neq\1$.

This semiring is {\it idempotent} if (3.6) is valid.  Idempotent semirings are
often called dioids. A semiring (maybe non-idempotent) is called {\it
  commutative}, if (3.7) is valid. Note that different versions of this
axiomatics are used, see, for example,
\cite{4}--\cite{9},\,\cite{14},\,\cite{20}--\cite{24}, \cite{27},\,\cite{28}
and some literature indicated in~\cite{53}.

{\bf Example 3.1.} The set $\Bbb R_+$ of all nonnegative real numbers endowed
with the usual addition and multiplication is a commutative (but not
idempotent) semiring. Of course, the field $\Bbb R$ of all real numbers is also
a commutative semiring.

{\bf Example 3.2.} $\Bbb R_{\max}$ and $\Bbb R_{\min}$ are isomorphic
commutative idempotent semirings.

{\bf Example 3.3.} $A=\Bbb R_+$ with the operations $\oplus=\max$ and
$\odot=\cdot$ (the usual multiplication); $\0=0$, $\1=1$. This idempotent
semiring is isomorphic to $\Bbb R_{\max}$ by the mapping $x\mapsto\ln(x)$.

{\bf Example 3.4.} $A=[a,b]=\{ x\in\Bbb R\vert\; a\leq x\leq b\}$ with the
operations $\oplus=\max$, $\odot=\min$ and the neutral elements $\0=a$ and
$\1=b$ (the cases $a=-\infty$, $b=+\infty$ are possible). $\bullet$

Semirings similar to these examples are the most close to the initial
``quantum'' object $\Bbb R_+$ and can be obtained by dequantization procedures.
However there are many important idempotent semirings which are unobtainable by
means of these procedures. Note that there exist important quantum mechanical
systems which cannot be obtained from classical systems by quantization (for
example, particles with spin and systems consisting of identical particles).
Thus the situation is natural enough for our analogy.

{\bf Example 3.5.} Let $\Mat_n(A)$ be the set of $n\times n$ matrices with
entries belonging to an idempotent semiring $A$.  This set forms a
noncommutative idempotent semiring with respect to matrix addition $\oplus$ and
matrix multiplication $\odot$, that is
\begin{displaymath}
  (X\oplus Y)_{ij} =X_{ij}\oplus Y_{ij} \qquad \mbox{and}\qquad (X\odot Y)_{ij}
  = \oplus_{k=1}^n X_{ik}\odot Y_{kj}.
\end{displaymath}
Of course, $(\0)_{ij}=\0\in A$, and $(\1)_{ij}=\0\in A$ if $i\neq j$, and
$(\1)_{ii}=\1\in A$.

{\bf Example 3.6.} $A=\{ 0,1\}$ with the operations $\oplus=\max$,
$\odot=\min$, $\0=0$, $\1=1$. This is the well-known {\it Boolean} semiring (or
Boolean algebra). $\bullet$

Note that every bounded distributive lattice is an idempotent semiring.

{\bf Example 3.7.} $A=\{\0,\1,a\}$, where $\{\0,\1\}$ is a Boolean semiring,
$\0\oplus a=a$, $\0\odot a=\0$, $\1\odot a=a$, $\1\oplus a=\1$, $a\oplus a=a$,
$a\odot a=a$. This example can be treated as a three-valued logic. $\bullet$

There are many finite idempotent semirings; a classification of commutative
idempotent semirings consisting of two, or three, or four elements is presented
in \cite{52}.

{\bf Example 3.8.} Let $A$ be the set of all compact convex subsets of $\Bbb
R^n$ (or of any closed convex cone in $\Bbb R^n$); this set is an idempotent
semiring with respect to the following operations:
\begin{eqnarray*}
  \alpha\oplus\beta&=&\mbox{ convex}\; \mbox{hull}\; \mbox{of } \; \alpha\;
  \mbox{ and }\;\beta;\\ \alpha\odot\beta&=&\{ a+b \mid a\in \alpha,
  b\in\beta\}
\end{eqnarray*}
for all $\alpha,\beta\in A$; $\0=\oslash$, $\1=\{0\}$.  This idempotent
semiring is used in mathematical economics and in the multicriterial
optimization problem (evolution of the so-called Paret\'o sets; see, for
example \cite{35}, \cite{74}).

{\bf Example 3.9.} If $A_1$ and $A_2$ are idempotent semirings, then
$A=A_1\times A_2$ is also an idempotent semiring with respect to the natural
component-wise operations of the direct product; in this case $(\0,\0)$ and
$(\1,\1)$ are the corresponding neutral elements.  A similar (and natural, see
\cite{52}) construction turns
$(A_1\backslash\{\0\})\times(A_2\backslash\{\0\})\cup\0$ into an idempotent
semiring. $\bullet$

Probably the first interesting and nontrivial idempotent semiring of all
languages over a finite alphabet was examined by S. Kleene \cite{73} in 1956.
This noncommutative semiring was used for applications to compiling and syntax
analysis, see also \cite{6}, \cite{7}.  There are many other interesting
examples of idempotent semirings (including the so-called ``tropical''
semirings, see, for example, \cite{47}, \cite{48}, \cite{60}, \cite{63},
\cite{64}) with applications to theoretical computer science (linguistic
problems, finite automata, discrete event systems and Petri nets, stochastic
systems, computational problems etc.), algebra (semigroups of matrices over
semirings), logic, optimization etc.; in particular, see also
\cite{5}--\cite{7},\,\cite{9},
\cite{11},\,\cite{12},\,\cite{15}--\cite{17},\,\cite{19}--\cite{24},\,\
\cite{26}--\cite{29},\,\cite{32},\,\cite{33},\,\cite{35},\,\
\cite{53},\,\cite{63}--\cite{66}.

There is a naturally defined {\it partial order} (i.e. partial ordering
relation) on any idempotent semiring (as well as on any idempotent semigroup);
by definition, $\alpha\preceq b$ if and only if $a\oplus b=b$. For this
relation the reflexivity is equivalent to the idempotency of the (generalized)
addition, whereas the transitivity and the antisymmetricity follow,
respectively, from the associativity and from the commutativity of this
operation.  This ordering relation on $\Bbb R_{\max}$ (as well as on semirings
described in the examples 3.3 and 3.4) coincides with the natural one but for
$\Bbb R_{\min}$ it is opposite to the natural ordering relation on the real
axis.

Every element $a$ of an idempotent semiring $A$ is ``nonnegative'': $\0\preceq
a$; indeed, $\0\oplus a=a$ because of (3.2). Similarly, for all $a,b,c\in A$ we
have $a\oplus c\preceq b\oplus c$, and $a\odot c\preceq b\odot c$ if $a\preceq
b$.

Using this standard partial order it is possible to define in the usual way the
notions of upper and lower bounds, bounded sets, $\sup M$ and $\inf N$ for
upper/lower bounded sets $M$ and $N$ etc.  On the basis of these concepts an
algebraic approach to the subject is developed, see, for example,
\cite{4}--\cite{9}, \cite{17}, \cite{19}--\cite{24}, \cite{27}, \cite{32},
\cite{33}, \cite{52}, \cite{53}.

An idempotent semiring can be a metric or topological space with natural
correlations between topological and algebraic properties.  For example, for
$\Bbb R_{\min}$ there is a natural metric $\rho(x, y)=\vert
e^{-x}-e^{-y}\vert$, and for the semiring from Example 3.4 it is convenient to
use the metric $\rho(x, y)=\vert\arctan x-\arctan y\vert$ \,if $a=-\infty$,
$b=+\infty$. The corresponding ``topological'' approach was developed, e.g. in
\cite{14}, \cite{15}, \cite{19}--\cite{24}, \cite{39}, \cite{42},
\cite{61}--\cite{63}, \cite{66}--\cite{70}.

\section{Semirings with special properties}

\quad It is convenient to treat some special classes of semirings for which
some additional conditions are fulfilled. Let us discuss some conditions of
this type.

Suppose $A$ is an arbitrary semiring. The so-called {\it cancellation
  condition} is fulfilled for $A$ if $b=c$ whenever $a\odot b=a\odot c$ and
$a\neq \0$.  If the multiplication in $A$ is invertible on $A\backslash\{\0\}$,
then $A$ is called {\it a semifield}. Of course, the cancellation condition is
fulfilled for all semifields. For example, $\Bbb R_{\max}$ is a semifield.
Idempotent semirings with the cancellation condition or with an idempotent
multiplication are especially interesting.

For arbitrary commutative idempotent semirings with the cancellation condition
the following version of Newton's binomial formula is valid:
\begin{equation}
  (a\oplus b)^n=a^n\oplus b^n,
\end{equation}
see \cite{32}, \cite{33}. However, this formula is valid also for semirings
from Example 3.4 which have no the cancellation condition. It is easily proved
(by induction) that for arbitrary commutative idempotent semirings this
binomial formula has the form
\begin{equation}
  (a\oplus b)^n=\bigoplus^n_{i=0}a^{n-i}\odot b^i.
\end{equation}

Suppose $A$ is an arbitrary idempotent semiring.  Applying (4.2) to the
semiring generated by elements $\1, a\in A$, we deduce the following formula:
\begin{equation}
  (\1\oplus a)^n=\1\oplus a\oplus a^2\oplus\cdots\oplus a^n.
\end{equation}

Now let $A$ be an arbitrary semiring (maybe non-idempotent) and suppose that
the following infinite sum
\begin{equation}
  a^*=\bigoplus^{\infty}_{i=0}a^i=\1\oplus a\oplus a^2\oplus\cdots\oplus a^n
  \cdots
\end{equation}
is well-defined for an element $a\in A$. For concrete semirings $a^*$ may be
defined, e.g. as $\sup\limits_n\{(\1+a)^n\}$ or
$\lim\limits_{n\to\infty}(\1\oplus a)^n$.  This important star operation
$a\mapsto a^*$ was introduced by S. Kleene \cite{73}; the element $a^*$ is
called a {\it closure} of~$a$.

It is natural to set $a^*=(\1-a)^{-1}$ if $A$ is a field and $a\neq \1$. It is
easy to prove that $a^*=\1$, if $A$ is an idempotent semiring and $a\preceq\1$.
For $\Bbb R_{\max}$ the closure $a^*$ is not defined if $\1\prec a$.
The situation can be corrected if we add an element $\infty$ such that
$a\oplus\infty = \infty$ for all $\in\Bbb R_{\max}$, $\0\odot\infty=\0$,
$a\odot\infty=\infty$ for all $a\neq \0$. For this new semiring
$\overline{\Bbb R}_{\max}=\Bbb R_{\max}\cup\{\infty\}$ we have
$a^*=\infty$ if $\1\prec a$, see e.g. \cite{18}, \cite{30}. For all
semirings described in the examples 3.4, 3.6 and 3.7 we have $a^*=\1$ for any
element~$a$.

An idempotent semiring $A$ is {\it algebraically closed} (with respect to the
operation $\odot$) if the equation $x^n=a$ (where $x^n=x\odot\cdots\odot x$)
has a solution $x\in A$ for any $a\in A$ and any positive integer $n$, see
\cite{32}, \cite{33}. It is remarkable and important that the semiring $\Bbb
R_{\max}$ is algebraically closed in this sense. However, the equation
$x^2\oplus\1=\0$ has no solutions.

\section{Correspondence principle}

\quad The analogy with Quantum Mechanics discussed in section 2 leads to the
following {\it correspondence principle} in idempotent calculus:

{\it There is a {\rm(}heuristic{\rm)} correspondence between important, useful
  and interesting constructions and results over the field of real (or complex)
  numbers {\rm(}or the semiring of all nonnegative numbers{\rm)} and similar
  constructions and results over idempotent semirings in the spirit of the
  correspondence principle in Quantum Mechanics}.

{\bf Example 5.1. Semimodules} (see e.g. \cite{4}--\cite{12}, \cite{17},
\cite{22}--\cite{24}, \cite{27}, \cite{32}, \cite{33}, \cite{50}--\cite{53}).
{\it A set} $V$ {\it is called a semimodule over a semiring} $A$ (or an
$A$-semimodule), if there is a commutative associative addition operation
$\oplus$ in $V$ with neutral element $\0$, and a multiplication $\odot$ of
elements from $V$ by elements of $A$ is defined, and the following properties
are fulfilled:
\begin{eqnarray*}
  &(\lambda\odot\mu)\odot v=\lambda\odot(\mu\odot v) &\qquad\mbox{ for
    all}\;\lambda, \mu\in A,\; v\in V;\\ &\lambda\odot(v_1\oplus
  v_2)=\lambda\odot v_1\oplus\lambda\odot v_2 &\qquad\mbox{ for
    all}\;\lambda\in A,\; v_1, v_2\in V;\\
&\0\odot v=\lambda\odot\0=\0
  &\qquad\mbox{ for all}\;\lambda\in A,\; v\in V.
\end{eqnarray*}
The addition $\oplus$ in $V$ is assumed to be idempotent if $A$ is an
idempotent semiring (i.e. $v\oplus v=v$ for all $v\in V$).  Then we assume that
\begin{displaymath}
  \sup_{\alpha}\{\lambda_{\alpha}\}\odot v=
  \sup_{\alpha}\{\lambda_{\alpha}\odot v\}, \qquad \mbox{ if}\; v\in V \;
  \mbox{and}\; \sup_{\alpha} \{\lambda_{\alpha}\}\in A.
\end{displaymath}
\indent Roughly speaking, semimodules are ``linear spaces'' over semirings. The
simplest $A$-semimodule is the direct sum (product) $A^n=\{(a_1, a_2, \cdots,
a_n)\; : \; a_j\in A\}$.  The set of all endomorphisms $A^n\to A^n$ coincides
with the semiring $\Mat_n(A)$ of all $A$-valued matrices (see Example 3.5 in
Section 3).

The theory of $A$-valued matrices is an analog of the well-known
O.~Perron--G.~Frobenius theory of nonnegative matrices, see e.g. \cite{75}. For
example, let $A$ be an algebraically closed commutative idempotent semiring
with the cancellation condition and the sequence $a^n\oplus b$ stabilizes for
any $a\preceq\1$ and $b\neq \0$, $a,b\in A$. Then for every endomorphism $K$ of
$A^n$ ($n\geq 1$) there exists a nontrivial subsemimodule $S\subset A^n$ (an
``eigenspace'') and $\lambda\in A$ (an ``eigenvalue'') such that
$Kv=\lambda\odot v$
for all $v\in S$; this element $\lambda$ is unique if $K$ is irreducible, see
\cite{32}, \cite{33}. In particular, this result is valid if $A=\Bbb R_{\max}$
(or $\Bbb R_{\min}$). Similar results can be proved for semimodules of bounded
functions and continuous functions, see \cite{32}, \cite{33}, \cite{22} --
\cite{24}.  $\bullet$

Idempotent analysis deals with functions taking values in idempotent semirings
and with the corresponding function spaces (semimodules). Let $X$ be a set and
$A$ an idempotent semiring. Let us denote by $B(X,A)$ the set of all bounded
mappings (functions) $X\to A$ (i.e. mappings with order-bounded images)
equipped with a natural structure of an $A$-semimodule. If $X$ is finite, $X=\{
x_1,\dots,x_n\}$, then $B(X,A)$ can be identified with the semimodule $A^n$
(see Example 5.1 above). Actually $B(X,A)$ is an idempotent semiring with
respect to the corresponding pointwise operations.

Let $A$ be a metric semiring; then there is the corresponding uniform metric on
$B(X,A)$. Suppose that $X$ is a topological space and then denote by $C(X,A)$
the subsemimodule of continuous functions in $B(X,A)$.

Suppose now that the space $X$ is locally compact and then denote by
$C_{\0}(X,A)$ the $A$-semimodule of continuous $A$-valued functions with
compact supports endowed with a natural topology (see \cite{19} -- \cite{24}
for details).

These spaces (and some other spaces of this type) are examples of
``idempotent'' function spaces. Many basic ideas, constructions and results can
be borrowed to idempotent analysis from usual analysis and functional analysis.

{\bf Example 5.2. Idempotent integration and measures}.  For the sake of
simplicity set $A=\Bbb R_{\max}$ and let $X$ be a locally compact space. An
idempotent analog of the usual integration can be defined by the formula
\begin{equation}
  \int_X^{\oplus}\varphi(x)\; dx=\sup_{x\in X}\varphi(x),
\end{equation}
if $\varphi$ is continuous or upper semicontinuous function on $X$. The set
function
\begin{equation}
  m_{\varphi}(B)=\sup_{x\in B}\varphi(x),
\end{equation}
where $B\subset X$ is called an $A$-{\it measure} on $X$ and $m_{\varphi}(\cup
B_{\alpha})=\bigoplus\limits_{\alpha}m_{\varphi}
(B_{\alpha})=\sup\limits_{\alpha}m_{\varphi}(B_{\alpha})$, so the function
(5.2) is completely additive. An idempotent integral with respect to this
$A$-measure is defined by the formula
\begin{equation}
  \int_X^{\oplus}\psi(x)\; dm_{\varphi}=
  \int_X^{\oplus}\psi(x)\odot\varphi(x)\; dx= \sup_{x\in
    X}\psi(x)\odot\varphi(x).
\end{equation}
It is obvious that this integration is ``linear'' over $A$ and it is easy to
see that (5.1) and (5.3) can be treated as limits of Riemann's and Lebesgue's
sums. Of course, if $\oplus=\min$ for the corresponding semiring $A$, then
(5.3) turns into the formula
\begin{equation}
  \int_X^{\oplus}\psi(x)\; dm_{\varphi}=
  \int_X^{\oplus}\psi(x)\odot\varphi(x)\; dx= \inf_{x\in
    X}\psi(x)\odot\varphi(x).
\end{equation}
In this case, $\odot$ may coincide e.g. with max, or the usual addition or
multiplication. See \cite{14}, \cite{15}, \cite{19} -- \cite{24} for details.
$\bullet$

Note that in (5.4) we mean $\inf$ (i.e. the greatest lower bound) with respect
to the usual ordering of numbers. But if $\oplus=\min$, then this order is
opposite to the standard partial order defined for any idempotent semiring (see
Section~3 above). It is clear that (5.3) and (5.4) coincide from this point of
view. In general case $A$--measure and idempotent integral can be defined by
(5.2) and (5.3), e.g. if the corresponding functions are bounded and $A$ is
{\it boundedly complete}, i.e. every bounded subset $B\subset A$ has the least
upper bound $\sup B.$

There is a natural analogy between idempotent and probability measures.  This
analogy leads to a parallelism between probability theory and stochastic
processes on the one hand, and optimization theory and decision processes on
the other hand. That is why it is possible to develop optimization theory at
the same level of generality as probability and stochastic processes theory. In
particular, the Markov causality principle corresponds to the Bellman
optimality principle; so the Bellman principle is an $\Bbb R_{\max}$--version
of the Chapman--Kolmogorov equation for Markov stochastic processes, see e.g.
\cite{43}--\cite{46},\,\cite{24},\,\cite{26},\,\cite{29},\,\cite{56},\,\cite{63}.
Applications to the filtering theory can be found in \cite{44},\,\cite{46}.

{\bf Example 5.3. Group idempotent (convolution) semirings.} Let $G$ be a
group, $A$ an idempotent semiring; assume that $A$ is boundedly complete. Then
the space $B(G, A)$ of all bounded functions $G\to A$ (see above) is an
idempotent semiring with respect to the following idempotent analog
$\circledast$ of convolution:\\ \hbox to\hsize{\hfill
  $\displaystyle{(\varphi\circledast\psi)(g)=
    \int^{\oplus}_G\varphi(x)\odot\psi(x^{-1}\cdot g)dx.} $\hfill
  (5.5)$\bullet$} \setcounter{equation}{5} \indent Of course, it is possible to
consider other ``function spaces'' instead of $B(G, A).$ In \cite{23},
\cite{24} semirings of this type are referred to as {\it convolution
  semirings.}

{\bf Example 5.4. Fourier--Legendre transform,} see \cite{14}, \cite{3},
\cite{19}--\cite{24}.  Let $A=\Bbb R_{\max}, G=\Bbb R^n$ and $G$ is treated as
a group. The usual Fourier--Laplace transform is defined by the formula
\begin{equation}
  \varphi(x)\mapsto\tilde\varphi(\xi)=\int_Ge^{i\xi\cdot x}\varphi (x)dx,
\end{equation}
where $e^{i\xi\cdot x}$ is a character of the group $G$, that is a solution of
the following functional equation:
$$ f(x+y)=f(x)f(y).$$ The corresponding idempotent analog ( for the case
$A=\Bbb R_{\max}$) has the form
$$ f(x+y)=f(x)\odot f(y)=f(x)+f(y),$$ so ``idempotent characters'' are linear
functionals $x\mapsto\xi\cdot x=\xi_1 x_1+\cdots +\xi_nx_n.$ Thus (5.6) turns
into the following transform:
\begin{equation}
  \varphi(x)\mapsto\tilde\varphi(\xi)=\int^{\oplus}_G\xi\!\cdot\! x\odot
  \varphi(x)dx=\sup_{x\in G}(\xi\!\cdot\! x+\varphi (x)).
\end{equation}
This is the famous {\it Legendre transform.} Thus this transform is an $\Bbb
R_{\max}$--version of the Fourier--Laplace transform.$\bullet$

Of course, this construction can be generalized to different classes of groups
and semirings. Transformations of this type convert the generalized convolution
to pointwise multiplication and possesses analogs of some important properties
of the usual Fourier transform. For the case of semirings of Pareto sets the
corresponding version of the Fourier transform reduces the multicriterial
optimization problem to a family of singlecriterial problems \cite{35}.

The examples 5.3 and 5.4 can be treated as fragments of an idempotent version
of the representation theory. In particular, idempotent representations of
groups can be examined as representations of the corresponding convolution
semirings (i.e. idempotent group semirings) in semimodules.

According to the correspondence principle, many important concepts, ideas and
results can be converted from usual functional analysis to idempotent analysis.
For example, idempotent scalar product can be defined by the formula:
\begin{equation}
  (\varphi, \psi)=\int_X^{\oplus}\varphi(x)\odot\psi(x)dx,
\end{equation}
where $\varphi, \psi$ are $A$--valued functions belonging to a certain
idempotent function space. There are many interesting spaces of this type
including $B(X, A)$, $C(X, A)$, $C_{\0}(X, A)$, analogs of the Sobolev spaces
and so on. There are analogs for the well-known theorems of Riesz, Hahn--Banach
and Banach--Steinhaus; it is possible to treat dual spaces and operators, an
idempotent version of the theory of distributions (generalized functions) etc.;
see \cite{19}--\cite{24}, \cite{34}, \cite{36}, \cite{39}, \cite{40}, \cite{76}
for details.

{\bf Example 5.5. Integral operators.} It is natural to construct idempotent
analogs {\it of integral operators} in the form
\begin{equation}
  K: \varphi(y)\mapsto (K\varphi)(x)=\int_Y^{\oplus}K(x, y)\odot\varphi(y)dy,
\end{equation}
where $\varphi(y)$ is an element of a space of functions defined on a set $Y$
and taking their values in an idempotent semiring $A$, $(K\varphi)(x)$ is an
$A$--valued function on a set $X$ and $K(x, y)$ is an $A$--valued function on
$X\times Y$. If $A=\Bbb R_{\max},$ then (5.9) turns into the formula
\begin{equation}
  (K\varphi)(x)=\sup_{y\in Y}\{K(x, y)+\varphi(y)\}.
\end{equation}
Formulas of this type are standard for optimization problems, see e.g.
\cite{77}. $\bullet$

It is easy to see that the operator defined by (5.9) is linear over $A$, i.e.
$K$ is an $A$--endomorphism of the corresponding semimodule (function space).
Actually every linear operator acting in an idempotent function space and
satisfying some natural continuity--type conditions can be presented in the
form (5.9). This is an analog of the well--known L.~Schwartz kernel theorem.
The topological version of this result in spaces of continuous functions was
established in \cite{78}, \cite{76}; see also \cite{23}, \cite{24}. The
algebraic version of the kernel theorem for the space of bounded functions see
in \cite{32}, \cite{33} and (in a final form) in \cite{52}.

\section{Superposition principle}

\quad In Quantum Mechanics the correspondence principle means that the
Schr\"o\-din\-ger equation (which is basic for the theory) is linear. Similarly
in idempotent calculus the correspondence principle means that some important
and basic problems and equations (e.g. optimization problems, the Bellman
equation and its generalizations, the Hamilton--Jacobi equation) nonlinear in
the usual sense can be treated as linear over appropriate idempotent semirings,
see \cite{1}--\cite{3}, \cite{19}--\cite{24}.

{\bf Example 6.1. Idempotent dequantization for the heat equation.} Let us
start with the heat equation
\begin{equation}
  \pdd {u}{t}=\frac {h}{2}\frac {{\partial}^2u}{\partial x^2},
\end{equation}
where $x\in\Bbb R,\, t>0,$ and $h$ is a positive parameter.

Consider the following change of variables:\\ \centerline{$u\mapsto w=-h\ln
  u;$} it converts (6.1) to the following (integrated) version of the Burgers
equation:
\begin{equation}
  \pdd wt+\frac 12\left(\pdd wx\right)^2- \frac h2\frac
  {{\partial}^2w}{\partial x^2}=0.
\end{equation}
This equation is nonlinear but it can be treated as linear over the following
generalized addition $\oplus$ and multiplication $\odot$ (borrowed from the
usual addition and multiplication by the change of variables):
\begin{equation}
  w_1\oplus w_2=-h\ln (e^{-w_1/h}+e^{-w_2/h}),
\end{equation}
\begin{equation}
  w_1\odot w_2=w_1+w_2.
\end{equation}
So if $w_1$ and $w_2$ are solutions for (6.2), then their linear combination
with respect to the operations (6.3) and (6.4) is also a solution for this
equation. For $h\to 0$ (6.2) turns into a special case of the Hamilton--Jacobi
equation:
\begin{equation}
  \pdd wt+\frac 12\left(\pdd wx\right)^2=0.
\end{equation}
This is the dequantization procedure described in Section 2 above. So it is
clear that (6.3) and (6.4) turn into addition $\oplus=\min$ and multiplication
$\odot=+$ in the idempotent semiring $\Bbb R_{\min}$ and the equation (6.5) is
linear over $\Bbb R_{\min}$; thus the set of solutions for (6.5) is an $\Bbb
R_{\min}$--semimodule. This example was the starting point for the well--known
Hopf method of vanishing viscosity.$\bullet$

In general case the Hamilton--Jacobi equation has the following form:
\begin{equation}
  \pdd {S(x, t)}t+H\left(\textstyle{\pdd Sx}, x, t\right)=0,
\end{equation}
where $H$ is a smooth function on $\Bbb R^{2n}\times [0, T]$. Consider the
Cauchy problem for (6.6): $S(x, 0)=S_0(x),\, 0\le t\le T,\, x\in\Bbb R^n$.
Denote by $U_t$ the resolving operator, i.e. the map that assigns to each given
$S_0(x)$ the solution $S(x, t)$ of this problem at the moment of time $t$. Then
the map $U_t$ for each $t$ is a linear (over $\Bbb R_{\min}$) integral operator
in the corresponding $\Bbb R_{\min}$--semimodule.

The situation is similar for the Cauchy problem for the homogeneous Bellman
equation
$$\textstyle{\pdd St+H\left(\pdd Sx\right)=0,\quad S\vert_{t=0}=S_0(x),}$$
where $H\colon\Bbb R^n\to\Bbb R$ is a convex (not strictly) first order
homogeneous function
$$H(p)=\sup_{(f, g)\in V}(f\!\cdot\! p+g),\, f\in\Bbb R^n,\, g\in\Bbb R,$$ and
$V$ is a compact set in $\Bbb R^{n+1}.$ See \cite{23}, \cite{24}, \cite{39},
\cite{76} for details.

It is well--known that discrete version of the Bellman equation can be treated
as linear over idempotent semirings. The so-called {\it generalized stationary}
(finite dimensional) {\it Bellman equation} has the form
\begin{equation}
  S=HS\oplus F,
\end{equation}
where $S, H, F$ are matrices with elements from an idempotent semiring $A$ and
the corresponding matrix operations are described in Example 3.5 above (for the
sake of simplicity we write $HS$ instead of $H\odot S$); the matrices $H$ and
$F$ are given (specified) and it is necessary to determine $S$ from the
equation.

The equation (6.7) has the following solution:
\begin{equation}
  S=H^*F,
\end{equation}
where $H^*$ is the closure of $H\in \Mat_n(A)$, see Section 4 and Example 3.5
above. Recall that
\begin{equation}
  H^*=\1\oplus H\oplus H^2\oplus\dots\oplus H^k\oplus\dots\ ,
\end{equation}
if the right-hand side of (6.9) is well-defined. In this case $H^*=\1\oplus
HH^*$, so $H^*F=F\oplus HH^*F$; thus (6.8) is a solution of (6.7). For example,
if the sequence $H^{(N)}=\sum_{k=0}^N H^k$ stabilizes (i.e. there exists $N_0$
such that $H^{(N)}=H^{(N_0)}$ for all $N\geq N_0$), then (6.9) is well-defined
and can be calculated by means of a finite set of operations (steps).

This consideration and a version of the Gauss elimination method for solving
(6.7) were presented by S.\,Kleene \cite{73} in the case of the semiring of all
languages over a finite alphabet. B.\,A.\,Carre \cite{4} used semirings to show
that many important problems for graphs can be formulated in a unified manner
and are reduced to solving systems of algebraic equations. For example,
Bellman's method of solving shortest path problems corresponds to a version of
the Jacobi method for solving (6.7), whereas Ford's algorithm corresponds to a
version of the Gauss-Seidel method. The further development of this subject see
in \cite{4}--\cite{18}, \cite{21}--\cite{24}, \cite{27}--\cite{31}, \cite{53},
\cite{65}.

Let $A$ be a semiring (maybe non-idempotent). For each square $n\times n$
matrix $H=(h_{ij})\in \Mat_n(A)$ there is a standard way to construct a
geometrical object called a {\it weighted directed graph}. This object consists
of a set $X$ of $n$ elements $x_1, x_2,\dots,x_n$ together with the subset
$\cyrg$ of all ordered pairs $(x_i,x_j)\in X\times X$ such that $h_{ij}\ne\0$
and the mapping $h\colon\cyrg\to A\backslash\{\0\}$ given by the correspondence
$(x_i,x_j)\mapsto h_{ij}$. The elements of $X$ are called {\it nodes\/}, and
the members of $\cyrg$ are called {\it arcs\/}; $h_{ij}$ are arc {\it
  weights\/}.

In other terms the quadruple $M(X,\cyrg,h,A)$ can be treated as a discrete
medium with the points $x_i$, the set $\cyrg$ of links and the so-called link
characteristics $h$. This concept is convenient for analysis of parallel
computations and for synthesis of computing media. Mathematical aspects of
these problems are examined in \cite{14}; the further development of the
subject is presented e.g. in \cite{15}, \cite{81}; see also \cite{23},
\cite{24}, \cite{27}, \cite{29}--\cite{31}, \cite{61}. For example, the
operating period evaluation problem for parallel algorithms and digital
circuits leads to shortest path problems for $M(X,\cyrg,h,A)$, where $A=\Bbb
R_{\max}$.

Recall that a sequence of nodes and arcs of the form
\begin{equation}
  p=(y_0,a_1,y_1,a_2,y_2,\dots,a_k,y_k),
\end{equation}
where $k\geq0$, $y_i$ are nodes of the graph, and $a_i$ are arcs satisfying
$a_i=(y_{i-1},y_i)$, is called a {\it path\/} (of order $k$) from the node
$y_0$ to the node $y_k$ in $M(X,\cyrg,h,A)$. The {\it weight\/} $h(p)$ {\it of
  the path\/} (6.10) is a product of the weights of its arcs:
\begin{equation}
  h(p)=h(a_1)\odot h(a_2)\odot\dots \odot h(a_k).
\end{equation}

The so-called {\it Algebraic Path Problem\/} is to find the following matrix
$D=(d_{ij})$:
\begin{equation}
  d_{ij}\stackrel{\rm def}{=}\oplus_p h(p),
\end{equation}
where $i,j=1,2,\dots,n$, and $p$ runs through all paths from $x_i$ to $x_j$.  A
solution to this problem does not need to exist (the set of weights in (6.12)
may be infinite). However, if there exists a closure $H^*$ of the matrix
$H=(h_{ij})$, then the matrix
\begin{equation}
  D=(d_{ij})=H^*=\1\oplus H\oplus H^2 \oplus \dots \oplus H^k\oplus \dots
\end{equation}
can be treated as a solution of this problem. Moreover, $H^k$ corresponds to
the value $\oplus_p h(p)$, where $p$ contains exactly $k$ arcs. For example,
$h^{(2)}_{ij}=\oplus_{k=1}^n h_{ik}\odot h_{kj}$ are elements (coefficients)
of $H^2$, and each coefficient $h^{(2)}_{ij}$ corresponds to $\oplus_p h(p)$,
where p runs through paths from $x_i$ to $x_j$ with exactly two arcs;
similarly, $H^3=H^2\odot H$, etc.

{\bf Example 6.2. The shortest path problem.} Let $A=\Bbb R_{\min}$, so
$h_{ij}$ are real numbers. In this case
$$ d_{ij}=\oplus_p h(p)= \min_p h(p),$$ where (6.11) has the form
$$h(p)=h(a_1)+h(a_2)+\dots+h(a_k).$$

{\bf Example 6.3. The relation closure problem.} Let $A$ be the Boolean
semiring (see Example 3.6 in Section 3 above). In this case $H$ corresponds to
a relation $R\subset X\times X$, $h_{ij}$ being $\1$ if and only if the
relation holds between $x_i$ and $x_j$. Then the transitive and reflective
closure $R^*$ of the relation $R$ corresponds to the matrix $D=H^*$.

{\bf Example 6.4. The maximal (minimal) width path problem.} Let $A$ be a
semiring $\Bbb R\cup\{-\infty\}\cup\{\infty\}$ with the operations
$\oplus=\max$ and $\odot=\min$ (see Example 3.4). Then
\begin{displaymath}
  d_{ij}=\oplus_p h(p)=\max_p h(p),
\end{displaymath}
where $h(p)=\min\{h(a_1),h(a_2),\dots,h(a_k)\}$. If $h(a_i)$ is the width (or
channel capacity) of $a_i$, then $h(p)$ is the possible width (or channel
capacity) of $p$.

{\bf Example 6.5. The matrix inversion problem.} Let $A$ be the field $\Bbb R$
of real numbers (which is not an idempotent semiring). In this case
$$D=H^*=\1+H+H^2\dots=(1-H)^{-1},$$ if the series $\sum_{k=0}^{\infty} H^k$
converges; if the matrix $\1-H$ is invertible, then $(\1-H)^{-1}$ can be
treated as a ``regularized'' sum of this series; here $H^0=\1$ is the identity
matrix.

{\bf Example 6.6. A simple dynamic programming problem.} Let $A=\Bbb R_{\max}$,
so $h_{ij}$ are real numbers. Let us consider $h_{ij}$ as a {\it profit\/} of
moving from $x_i$ to $x_j$, and suppose $f_i$ is a {\it terminal prize\/} for
the node $x_i$ ($f_i\in\Bbb R$). Assume that $p$ is a path of the form (6.10)
and $y_0=x_i$. Let $M$ be a {\it total profit\/} for $p$, that is
$$M=h(a_1)+h(a_2)+\dots+h(a_k)+f(y_k).$$ It is easy to see that $\max
M=(H^kf)_i$, where $f$ is a vector $\{f_i\}$, $H,H^k\in \Mat_n(A)$. So, the
maximal value of the total profit for $k$ steps is $(H^kf)_i$. It is clear that
the maximal value of the total profit for paths of arbitrary order is $\max
M=(H^*f)_i$.$\bullet$

See many other examples and details (including semiring versions of linear
programming) in \cite{4}--\cite{17}, \cite{21}--\cite{24}, \cite{27},
\cite{30}, \cite{73}, \cite{79}, \cite{80}, \cite{82}. The book \cite{27}
of F.~L.~Baccelli, G.~Cohen, G.~J.~Olsder and J.-P.~Quadrat is
particularly useful.

\section{Correspondence principle for algorithms}

\quad Of course, the correspondence principle is valid for algorithms (as well
as for their software and hardware implementations). Thus:

{\it If we have an important and interesting numerical algorithm, then we have
  a good chance that its semiring analogs are important and interesting as
  well.}

In particular, according to the superposition principle, analogs of linear
algebra algorithms are especially important. Note that numerical algorithms for
standard infinite-dimensional linear problems over semirings (i.e. for problems
related to integration, integral operators and transformations, the
Hamilton--Jacobi and generalized Bellman equations) deal with the corresponding
finite-dimensional (or finite) ``linear approximations''.  Nonlinear algorithms
often can be approximated by linear ones. Recall that usually different natural
algorithms for the same optimization problem correspond to different standard
methods for solving systems of linear equations (like Gauss elimination method,
iterative methods etc.).

It is well-known that algorithms of linear algebra are convenient for parallel
computations (see, e.g. \cite{81}--\cite{84}); so, their idempotent analogs
accept a parallelization. This is a regular way to use parallel computations
for many problems including basic optimization problems.

Algorithms for the ``scalar'' (inner) product of two vectors, for matrix
addition and multiplication do not depend on concrete semirings. Algorithms to
construct the closure $H^*$ of an ``idempotent'' matrix $H$ can be derived from
standard methods for calculating $(\1-H)^{-1}$. For the Gauss--Jordan
elimination method (via LU-decomposition) this trick was used in \cite{30}, and
the corresponding algorithm is universal and can be applied both to the general
algebraic path problem and to computing the inverse of a real (or complex)
matrix $(\1-H)$. Computation of $H^{-1}$ can be derived from this universal
algorithm with some obvious cosmetic transformations.

Note that numerical algorithms are combinations of basic operations. Usually
these basic operations deal with ``numbers''. Actually these ``numbers'' are
thought as members of some numerical {\it domains\/} (real numbers, integers,
and so on). But every computer calculation deals with concrete {\it models\/}
(computer representations) of these numerical domains. For example, real
numbers can be represented as ordinary floating point numbers, or as double
precision floating point numbers, or as rational numbers etc. Differences
between mathematical objects and their computer models lead to calculation
errors. That is another reason to use universal algorithms which do not depend
on a concrete semiring and its concrete computer model. Of course, one
algorithm may be more universal than another algorithm of the same type. For
example, numerical integration algorithms based on the Gauss--Jacobi quadrature
formulas actually depend on computer models because they use finite precision
constants.  On the contrary, the rectangular formula and the trapezoid rule do
not depend on models and in principle can be used even in the case of
idempotent integration.

\section{Correspondence principle for hardware\newline design}

\quad A systematic application of the correspondence principle to computer
calculations leads to a unifying approach to software and hardware design.

The most important and standard numerical algorithms have many hardware
realizations in the form of technical devices or special processors.  {\it
  These devices often can be used as prototypes for new hardware units
  generated by substitution of the usual arithmetic operations for its semiring
  analogs and by addition tools for performing neutral elements\/} $\0$ {\it
  and\/} $\1$ (the latter usually is not difficult). Of course the case of
numerical semirings consisting of real numbers (maybe except neutral elements)
is the most simple and natural. Semirings of this type are presented in the
examples 3.1--3.4. Semirings from the examples 3.6 and 3.7 can also be treated
as numerical semirings. Note that for semifields (including $\Bbb R_{\max}$ and
$\Bbb R_{\min}$) the operation of division is also defined.

Good and efficient technical ideas and decisions can be transposed from
prototypes into new hardware units. Thus the correspondence principle generates
a regular heuristic method for hardware design. Note that to get a patent it is
necessary to present the so-called ``invention formula'', that is to indicate a
prototype for the suggested device and the difference between these devices. A
survey of patents from the correspondence principle point of view is presented
in [82].

Consider (as a typical example) the most popular and important algorithm of
computing the scalar product of two vectors:
\begin{equation}
  (x,y)=x_1y_1+x_2y_2+\dots+x_ny_n.
\end{equation}
The universal version of (8.1) for any semiring $A$ is obvious:
\begin{equation}
  (x,y)=(x_1\odot y_1)\oplus (x_2\odot y_2)\oplus \dots\oplus(x_n\odot y_n).
\end{equation}
In the case $A=\Bbb R_{\max}$ this formula turns into the following one:
\begin{equation}
  (x,y)=\max\{ x_1+y_1, x_2+y_2, \dots, x_n+y_n\}.
\end{equation}

This calculation is standard for many optimization algorithms (see Section 6),
so it is useful to construct a hardware unit for computing (8.3).  There are
many different devices (and patents) for computing (8.1) and every such device
can be used as a prototype to construct a new device for computing (8.3) and
even (8.2). Many processors for matrix multiplication and for other algorithms
of linear algebra are based on computing scalar products and on the
corresponding ``elementary'' devices respectively, etc.

There are some methods to make these new devices more universal than their
prototypes. There is a modest collection of possible operations for standard
numerical semirings: $\max$, $\min$, and the usual arithmetic operations. So,
it is easy to construct programmable hardware processors with variable basic
operations. Using modern technologies it is possible to construct cheap
special-purpose multi-processor chips implementing examined algorithms. The
so-called systolic processors are especially convenient for this purpose. A
systolic array is a ``homogeneous'' computing medium consisting of elementary
processors, where the general scheme and processor connections are simple and
regular. Every elementary processor pumps data in and out performing elementary
operations in a such way that the corresponding data flow is kept up in the
computing medium; there is an analogy with the blood circulation and this is a
reason for the term ``systolic'', see e.g.  \cite{83}, \cite{84}.

Concrete systolic processors for the general algebraic path problem are
presented in \cite{30}, \cite{31}. In particular, there is a systolic array of
$n(n+1)$ elementary processors which performs computations of the Gauss--Jordan
elimination algorithm and can solve the algebraic path problem within $5n-2$
time steps. Of course, hardware implementations for important and popular basic
algorithms increase the speed of data processing.

\section{Correspondence principle for software\newline design}

\quad Software implementations for universal semiring algorithms are not so
efficient as hardware ones (with respect to the computation speed) but are much
more flexible.  Program modules can deal with abstract (and variable)
operations and data types. Concrete values for these operations and data types
can be defined by input data types. In this case concrete operations
and data types are defined by means of additional program modules. For programs
written in this manner it is convenient to use a special techniques of the
so-called object-oriented design, see e.g. \cite{71}.
Fortunately, powerful tools supporting the object-oriented software
design have recently appeared including compilers for real and convenient
programming languages (e.g.  $C^{++}$).

There is a project to obtain an implementation of the correspondence principle
approach to scientific calculations in the form of a powerful software system
based on a unifying collection of universal algorithms.  This approach ensures
a working time reduction for programmers and users because of software
unification. The arbitrary necessary accuracy and safety of numerical
calculations can be ensured as well \cite{72}.

The system contains several levels (including the programmer and user levels)
and many modules. Roughly speaking it is divided into three parts.  The first
part contains modules that implement finite representations of basic
mathematical objects (arbitrary precision real and complex numbers, finite
precision rational numbers, $p$-adic numbers, interval numbers, fuzzy numbers,
basic semirings and rings etc.). The second part implements universal
calculation algorithms (linear algebra, idempotent and usual analysis,
optimization and optimal control, differential equations and so on). The third
part contains modules implementing model dependent algorithms (e.g. graphics,
Gauss--Jacobi type numerical integration, efficient approximation algorithms).
The modules can be used in user programs written in $C^{++}$. See \cite{72} for
some details.

E-mail: litvinov@islc.msk.su

glitvinov@mail.ru
\end{document}